\documentclass[a4paper,12pt]{article}
\usepackage[latin1]{inputenc}
\usepackage{amsmath,amsthm,amsfonts}
\usepackage{amssymb,amstext,amsgen}
\usepackage{amsbsy,amsopn}
\usepackage[pagebackref]{hyperref}
\usepackage{mathrsfs}
\usepackage{amscd}
{} 
{Tensor Distributions
 \hfill
 \qquad M. Grosser\qquad}
\begin{document}
\newcommand{\TO}{\mathop{\mathrm{TO}}}
\title{Tensor Valued Colombeau Functions\\ on Manifolds\footnote{This work was
partially supported by projects P20525 and Y237 of the 
Austrian Science Fund.}}
\author{M. Grosser\footnote{Faculty of Mathematics,
University of Vienna, Nordbergstra\ss e 15, A-1090, Austria.
e-mail: michael.grosser@univie.ac.at}}
\date{}
\maketitle
\begin{abstract}Extending the construction of the algebra $\hat{\mathcal G}(M)$
of
scalar valued Colombeau functions on a smooth manifold $M$ (cf.\
[4]),
we present a suitable basic space for eventually obtaining
tensor valued generalized functions on $M$, via the usual quotient construction.
This basic space canonically contains the tensor valued distributions and
permits a natural extension of the classical Lie derivative. Its members are
smooth functions depending---via a third slot---on so-called transport
operators, in addition to slots one (smooth $n$-forms on $M$) and two (points of
$M$) from the scalar case.
\\[4mm] {\it AMS Mathematics  Subject Classification (2000)}:
Primary 46F30; Secondary 46T30, 53A45.
\\[2mm] {\it Key words and phrases:} Colombeau functions, tensors,
manifolds.
\end{abstract}

\vskip24pt
In the following, we generalize the construction of the full Colombeau
algebra $\hat{\mathcal G}(M)$ (see
[4])
%\cite{steingf07})
to the
tensor valued case. Let $M$ denote an (orientable) smooth paracompact Hausdorff
manifold of dimension $n$; always let $p\in M$, $f\in{\mathcal C}^\infty(M)$,
$u\in{\mathcal
D}'(M)$, $\omega\in\hat{\mathcal A}_0(M)(\subseteq\Omega^n_c(M))$;
$R\in\hat{\mathcal E}(M)$ (notation as in
[4]).
%\cite{steingf07}).
Recall
the scalar case setting:
% lined out in
%([4]):
%%\cite{steingf07}):

\vskip-5pt
$$\begin{array}{ccccc}
\mathrm{Smooth\ functions}:\ f\in {\mathcal C}^\infty(M)
     &&\kern-30pt\mathrm{Distributions}:\ u\in{\mathcal D}'(M)\\[4pt]
\fbox{$f(p)$}&&\ \kern-30pt \fbox{$u(\omega)=\langle u,\omega\rangle$}\\[7pt]
\mathrm{SLOT\ 1} && \kern-30pt\mathrm{SLOT\ 2}\\[9pt]
& \kern-30pt\mathrm{Colombeau}&\\
& \kern-30pt\mathrm{generalized\ functions}:&\\
& \kern-30pt R\in \hat{\mathcal E}(M)&\\[4pt]
& \kern-30pt\fbox{$R(\omega,p)$}&\\[7pt]
& \kern-30pt\mathrm{SLOTS\ 1,2} &
\end{array}$$
%

%\pagebreak
Embedding smooth functions $f$ by $\sigma$ resp.\ distributions $u$ by $\iota$
into
$\hat{\mathcal E}(M)$ is effected by using slot 1 resp.\ slot 2, by means of
the formulas (well-known from
[2]
%\cite{gksv}
resp.\
[4])
%\cite{steingf07}
\begin{eqnarray}
\label{sigma} (\sigma f)(\omega,p)&:=&f(p)\\
\label{iota}  (\iota u) (\omega,p)&:=&u(\omega)=\langle u,\omega\rangle.
\end{eqnarray}
%For the time being, i.e.\ before passing to quotients of moderate by negligible
%elements, we ignore the fact that ${\mathcal C}^\infty(M)$ in fact is a
%subset of ${\mathcal D}'(M)$.

Starting from $\hat{\mathcal E}(M)$, the Colombeau algebra $\hat{\mathcal G}(M)$
is then constructed
%as some quotient of some subalgebra of it,
by passing to quotients of moderate by negligible
elements, as usual in
Colombeau theory.
On the level of quotients resp.\ classes, $\sigma$ and $\iota$ become
equal on ${\mathcal C}^\infty(M)$. However, we do not actually perform this last
step of the construction at the moment, the
question of appropriate basic spaces being our main focus.

For a long period the guiding intuitive idea of the authors of
[2]
%\cite{gksv}
towards obtaining a suitable basic space for tensor valued
generalized functions on $M$ had been the following:

\vskip6pt
for {\em scalars} on $M$ use $\hat{\mathcal E}(M)$, the candidate 
which had proven successful in
[2];
%\cite{gksv};

\vskip3pt
for {\em tensors} on $M$ perform an appropriate ``afterward'' tensorial
construction based on the ready-made space $\hat{\mathcal E}(M)$.

\vskip6pt
All efforts along these lines essentially led to some version of
``coordinate-wise embedding'' $\iota^r_s$ of distributional tensor fields of
type $(r,s)$ ($r$ contravariant, $s$ covariant indices).
%; see below for details.
This way of proceeding, however, is ultimately barred due to a consequence of
the
famous Schwartz type impossibility result: Viewing $\iota$ as a map embedding
${\mathcal D}'(M)$ into $\hat{\mathcal G}(M)$ as in
[2],
%\cite{gksv},
we have, in
general,
\begin{equation}\label{different}
\iota(fu)\neq\iota(f)\cdot\iota(u)\qquad\qquad(f\in{\mathcal C}^\infty(M),\
u\in{\mathcal D}'(M)),
\end{equation}
that is, $\iota$ is not ${\mathcal C}^\infty(M)$-linear.

\pagebreak
To get an impression of what a tensorial construction
as just indicated should look like and in which way
the above Schwartz type result poses an unsurmountable obstacle to the approach
of coordinate-wise embedding we review the situation for tensorial distributions
(of type $(r,s)$, say) on $M$. To this end, denote by ${\mathrm T}^r_sM$ the
bundle of $(r,s)$-tensors over $M$ and by ${\mathcal T}^r_s(M)$ the linear space
of smooth sections of ${\mathrm T}^r_sM$, i.e.\ of smooth tensor fields of type
$(r,s)$ on $M$. The linear space ${{\mathcal D}'}^r_s(M)$ of tensorial
distributions of type $(r,s)$ on $M$ can be defined in several equivalent ways;
for our present purpose, we prefer
$${{\mathcal D}'}^r_s(M) := 
          ({\mathcal T}^s_r(M)\otimes_{{\mathcal C}^\infty(M)}\Omega^n_c(M))'
$$
(compare section
3.1.3 of
[1]
%\cite{gkos}
where---due to not assuming orientability of
$M$---densities on $M$ take the place of $n$-forms, yielding a slightly more
general setting).

Now it is a fundamental result that tensorial distributions can be viewed as
tensor fields with (scalar) distributional coefficients
([1],
%(\cite{gkos},
3.1.15), i.e.,
\begin{equation}\label{distrtensor}
{{\mathcal D}'}^r_s(M)\cong{\mathcal D}'(M)
                      \otimes_{{\mathcal C}^\infty(M)}{\mathcal T}^r_s(M).
\end{equation}
A formula completely analogous to (\ref{distrtensor}) is valid (though trivial)
on the level of smooth objects:
\begin{equation}\label{smoothtensor}
{\mathcal T}^r_s(M)\cong{\mathcal C}^\infty(M)
                      \otimes_{{\mathcal C}^\infty(M)}{\mathcal T}^r_s(M).
\end{equation}
(\ref{distrtensor}) and (\ref{smoothtensor}) are interlaced by natural
isomorphisms:
Denoting the embedding of smooth regular objects into distributional ones as

%\pagebreak
\begin{eqnarray*}
\rho&:&{\mathcal C}^\infty(M)\to{\mathcal D}'(M)\\
\rho^r_s&:&{\mathcal T}^r_s(M)\to{{\mathcal D}'}^r_s(M)
\end{eqnarray*}
we obtain the following commutative ``TD-diagram'':

\vskip9pt
\[
\begin{CD}\label{TDdiagram}
{\mathcal T}^r_s(M) @>\cong>> {\mathcal C}^\infty(M)\otimes_{{\mathcal
     C}^\infty(M)}
             {\mathcal T}^r_s(M) \\
@V\rho^r_sVV @VV\rho\otimes
             \mathrm{id}V\\
{{\mathcal D}'}^r_s(M) @>\cong>>{\mathcal D}'(M)
\otimes_{{\mathcal
      C}^\infty(M)}{\mathcal T}^r_s(M) 
\end{CD}
\]

\vskip9pt\noindent
This certainly encourages us to try the definition
\begin{equation}\label{generaltensor}
\hat{\mathcal G}^r_s(M):=\hat{\mathcal G}(M)
                      \otimes_{{\mathcal C}^\infty(M)}{\mathcal T}^s_r(M)
\end{equation}

\pagebreak\noindent
yielding the reassuring ``TG-diagram''

\vskip9pt
\[
\begin{CD}\label{TGdiagram}
{\mathcal T}^r_s(M) @>\cong>> {\mathcal C}^\infty(M)\otimes_{{\mathcal
     C}^\infty(M)}
             {\mathcal T}^r_s(M) \\
@V\sigma^r_s VV @VV
          \sigma\otimes\mathrm{id}V\\
\hat{\mathcal G}^r_s(M) @>\cong>>\hat{\mathcal G}(M)
\otimes_{{\mathcal
      C}^\infty(M)}{\mathcal T}^r_s(M) 
\end{CD}
\]

\vskip9pt\noindent
Combining the TD- and the TG-diagrams into one (and omitting the
${\mathcal C}^\infty(M)$-subscript at the $\otimes$ sign, as well as all
occurrences of ``$(M)$'') results in

\vskip9pt
\begin{center}
\setlength{\unitlength}{1500sp}%
\begingroup\makeatletter\ifx\SetFigFont\undefined%
\gdef\SetFigFont#1#2#3#4#5{%
  \reset@font\fontsize{#1}{#2pt}%
  \fontfamily{#3}\fontseries{#4}\fontshape{#5}%
  \selectfont}%
\fi\endgroup%
\begin{picture}(8550,6747)(1201,-6973)
\thinlines
\put(2401,-961){\vector( 1, 0){6000}}
\put(1801,-1561){\vector( 0,-1){4800}}
\put(9001,-1486){\vector( 0,-1){4875}}
\put(2401,-6961){\vector( 1, 0){6000}}
\put(2401,-1561){\vector( 2,-3){1200}}
\put(8401,-1561){\vector(-2,-3){1165.385}}
\put(4201,-3961){\vector( 1, 0){2400}}
\put(3601,-4561){\vector(-2,-3){1200}}
\put(7201,-4561){\vector( 2,-3){1200}}
\put(1501,-961){\makebox(0,0)[lb]{\smash{\SetFigFont{12}{14.4}{\rmdefault}
{\mddefault}{\updefault}${\mathcal T}^r_s$}}}
\put(8901,-961){\makebox(0,0)[lb]{\smash{\SetFigFont{12}{14.4}{\rmdefault}
{\mddefault}{\updefault}${\mathcal C}^\infty\otimes{\mathcal T}^r_s$}}}
\put(5001,-711){\makebox(0,0)[lb]{\smash{\SetFigFont{12}{14.4}{\rmdefault}
{\mddefault}{\updefault}$\cong$}}}
\put(5001,-3661){\makebox(0,0)[lb]{\smash{\SetFigFont{12}{14.4}{\rmdefault}
{\mddefault}{\updefault}$\cong$}}}
\put(5001,-6661){\makebox(0,0)[lb]{\smash{\SetFigFont{12}{14.4}{\rmdefault}
{\mddefault}{\updefault}$\cong$}}}
\put(1651,-6961){\makebox(0,0)[lb]{\smash{\SetFigFont{12}{14.4}{\rmdefault}
{\mddefault}{\updefault}$\hat{\mathcal G}^r_s$}}}
\put(8901,-6961){\makebox(0,0)[lb]{\smash{\SetFigFont{12}{14.4}{\rmdefault}
{\mddefault}{\updefault}$\hat{\mathcal G}\otimes{\mathcal T}^r_s$}}}
%mg:
\put(3363,-2500){\makebox(0,0)[lb]{\smash{\SetFigFont{12}{14.4}{\rmdefault}
{\mddefault}{\updefault}\hspace*{-.5mm}$\rho^r_s$}}}
\put(7189,-2500){\makebox(0,0)[lb]{\smash{\SetFigFont{12}{14.4}{\rmdefault}
{\mddefault}{\updefault}\hspace*{-7.5mm}%\iota\otimes\mathrm{id}
$\rho\otimes\mathrm {id}$}}}
%mg ende
%mike:
\put(3363,-5561){\makebox(0,0)[lb]{\smash{\SetFigFont{12}{14.4}{\rmdefault}
{\mddefault}{\updefault}\hspace*{-.5mm}?\,\fbox{1}\,?}}}
\put(7189,-5561){\makebox(0,0)[lb]{\smash{\SetFigFont{12}{14.4}{\rmdefault}
{\mddefault}{\updefault}$\hspace*{-7.5mm}%\iota\otimes\mathrm{id}
?\,\fbox{2}\,?$}}}
%mike ende
\put(1201,-4100){\makebox(0,0)[lb]{\smash{\SetFigFont{12}{14.4}{\rmdefault}
{\mddefault}{\updefault}\hspace*{-1mm}$\sigma^r_s$}}}
\put(9351,-4100){\makebox(0,0)[lb]{\smash{\SetFigFont{12}{14.4}{\rmdefault}
{\mddefault}{\updefault}$\sigma\otimes\mathrm{id}$}}}
\put(3401,-4100){\makebox(0,0)[lb]{\smash{\SetFigFont{12}{14.4}{\rmdefault}
{\mddefault}{\updefault}${\mathcal D}'^r_s$}}}
\put(7101,-4100){\makebox(0,0)[lb]{\smash{\SetFigFont{12}{14.4}{\rmdefault}
{\mddefault}{\updefault}$\hspace*{-2mm}{\mathcal D}'\otimes{\mathcal T}^r_s$}}}
\end{picture}
\end{center}

\vskip9pt\noindent
where the arrows denoted by \fbox{1} resp. \fbox{2} still are waiting to be
defined---the former providing the desired embedding of tensor distributions
into generalized tensors. Now, \fbox{1} certainly would have to be induced by
\fbox{2}\,,
and for the latter, due to $\sigma = \iota\circ \rho$, the only sensible choice
is
$\iota\otimes\mathrm{id}$. However, we have to remember that our $\otimes$ signs
actually read $\otimes_{{\mathcal C}^\infty(M)}$. Therefore, we have to check
carefully whether mappings giving rise to a commutative ``DG-diagram''

\vskip9pt
\[
\begin{CD}\label{DGdiagram}
{{\mathcal D}'}^r_s(M) @>\cong>> {\mathcal D}'(M)\otimes_{{\mathcal
     C}^\infty(M)}
             {\mathcal T}^r_s(M) \\
@V\iota^r_s\mathrm{\ ??} VV @VV
          \iota\otimes\mathrm{id\ ??}V\\
\hat{\mathcal G}^r_s(M) @>\cong>>\hat{\mathcal G}(M)
\otimes_{{\mathcal
      C}^\infty(M)}{\mathcal T}^r_s(M) 
\end{CD}
\]

\vskip9pt\noindent
actually exist. Unfortunately, the answer is no! To be sure, on the level
of vector space tensor products,
$$\iota\otimes\mathrm{id}:
   {\mathcal D}'(M)\otimes{\mathcal T}^r_s(M)\to
   \hat{\mathcal G}(M)\otimes {\mathcal T}^r_s(M)$$
is well-defined. Yet it does not induce a corresponding map on the level of 
${\mathcal C}^\infty(M)$-module tensor products (which would be what we
actually need) since it is not balanced into
$\hat{\mathcal G}(M)
     \otimes_{{\mathcal C}^\infty(M)} {\mathcal T}^r_s(M)$
by the Schwartz type theorem:
$$(\iota\otimes\mathrm{id})((f\cdot u)\otimes t)=
     \iota(f\cdot u)\otimes t$$
is different in general (cf.\ (\ref{different})) from
$$(\iota\otimes\mathrm{id})(u\otimes (f\cdot t))=
     \iota(u)\otimes (f\cdot t)=
     (\iota(u)\cdot f)\otimes t=
     (\sigma(f)\cdot\iota(u))\otimes t=
     (\iota(f)\cdot\iota(u))\otimes t
     .$$

It is instructive to take a look at the coordinate version of the preceding
(geometrically phrased) impossibility result. As we will show, the attempt to
%employ something like
build upon $\iota\otimes\mathrm{id}$ is reflected by trying to embed
tensor fields coordinate-wise. Again we will arrive at a
contradiction, demonstrating that
coordinate-wise embedding has to
be abandoned completely when spaces of tensor valued Colombeau
functions---allowing for a canonical
embedding of distributions---are to be constructed.

For localizing, assume that $M$ can be described by a single chart. Then 
${\mathcal T}^r_s(M)$ has a ${\mathcal C}^\infty(M)$-basis consisting of
smooth tensor fields, say, $e_1,\dots,e_m\kern-.3pt\in\kern-.3pt{\mathcal
T}^r_s(M)$ with
$m=n^{r+s}$. By
(\ref{generaltensor}), every $u\in{{\mathcal D}'}^r_s(M)$ can be written as
$u=u^i\otimes e_i$ (using summation convention) with $u^i\in{\mathcal D}'(M)$.
The geometrical requirement of $\iota\otimes\mathrm{id}$ being well-defined on
the level of module tensor products corresponds to 
$(\iota\otimes\mathrm{id})(u)$ being independent of basis representation of
$u$. Thus consider a change of basis given by $e_i=a^j_i\hat e_j$, with $a^j_i$
smooth. Then $u=\hat u^j\otimes \hat e_j$ with $\hat u^j = a^j_i u^i$. Applying 
$\iota\otimes\mathrm{id}$ to either representation of $u$, we obtain
$$(\iota\otimes\mathrm{id})(u^i\otimes e_i)\kern-3pt=\kern-3pt
     \iota (u^i)\otimes (a^j_i\hat e_j)\kern-3pt=\kern-3pt
     (\iota (u^i)a^j_i)\otimes \hat e_j\kern-3pt=\kern-3pt
     (\sigma(a^j_i)\iota(u^i))\otimes \hat e_j\kern-3pt=\kern-3pt
     (\iota(a^j_i)\iota(u^i))\otimes \hat e_j$$
resp.\
$$(\iota\otimes\mathrm{id})(\hat u^j\otimes \hat e_j)=
     \iota(a^j_i u^i)\otimes \hat e_j$$
which are different in general due to
$\iota(a^j_i)\iota(u^i)\neq \iota(a^j_i u^i)$
(cf.\ (\ref{different})).
It should be clear now that relying on coordinate-wise embedding is betting on
the wrong horse.

%\vskip1cm
%{\tt IMPOSSIBILITY RESULT}
%\vskip1cm
To circumvent this Schwartz type obstacle, the following alternative approach
(due mainly to J. A. Vickers and J. P. Wilson, cf.\
[5])
%\cite{vw})
turned out to be
successful eventually: Introduce, in addition to slots 1 and 2, some slot 3
``inside'' of $R$, i.e. intervene ``before'' $R$ actually acts by assigning
some tensor to its argument(s).

From now on, let us write ``$t$'' (for ``tensor'') rather than ``$R$''. Thus
the new idea directs us to replace $R(\omega,p)$ by
$$t(\omega,p,A)$$
($A$ having been fed into slot 3) in a way that $t$ becomes a member of some
space $\hat{\mathcal E}^r_s(M)$ of (smooth) tensor valued functions, to be
defined appropriately. This latter space then will serve as the basic space for
tensors of type $(r,s)$, consisting of functions having three slots as above.

Observe that this strategy includes ``redefining'' also the scalar case,
in a way that the ''old'' 2-slot version from
[2] resp.\ [4]
%\cite{gksv}
has to be upgraded
to the ``new'' 3-slot version. So, strictly speaking the algebra
$\hat{\mathcal G}(M)$ of
(scalar) valued generalized functions discussed in [2] resp.\ [4] in fact
differs (by the absence/presence of slot 3) from the algebra
$\hat{\mathcal G}^0_0(M)$
introduced (as the special case $r=s=0$ of $\hat{\mathcal G}^r_s(M)$) at
the end of this article.

Now let us explain and motivate which kind of objects we should expect to feed
into slot 3. As to $\omega$ and $p$, we take $\omega\in\hat{\mathcal A}_0(M)$
resp.\
$p\in M$, as we did previously for $R\in\hat{\mathcal E}(M)$. $A$, on the other
hand, has to be taken as a member of $\Gamma_c(\mathop{\mathrm{TO}}(M,M))$, the
latter denoting the
space of
compactly supported
smooth sections of the bundle $\mathop{\mathrm{TO}}(M,M)$ of
``transport operators''
over $M\times M$. More explicitly, $A$ is a compactly supported smooth map
$$A:M\times M\to\bigsqcup_{(p,q)\in M\times M}
{\mathrm L}({\mathrm T}_p M,{\mathrm T}_q M)$$
where ${\mathrm L}({\mathrm T}_p M,{\mathrm T}_q M)$ denotes the space of all
linear maps from the tangent space at $p$ to $M$ into the tangent space at $q$
of $M$, and the disjoint union above carries the bundle structure suggested
by the obvious local coordinate respresentations. Thus we have, for $p,q\in M$,
$$A(p,q):{\mathrm T}_p M\to{\mathrm T}_q M\qquad\mathrm{(linear)}$$
where $A(p,q)$ smoothly depends on $p$ and $q$.

%Denoting the bundle of tensors of type $(r,s)$ over $M$ by ${\mathrm T}^r_s M$,
The new basic space
$\hat{\mathcal E}^r_s(M)$ will
%then
be defined as a certain subspace (to be
specified later) of
$${\mathcal C}^\infty(\hat{\mathcal A}_0(M)\times M\times
\Gamma_c(\mathop{\mathrm{TO}}(M,M)),{\mathrm T}^r_s M),$$
or, with $\hat{\mathcal B}(M):=\Gamma_c(\mathop{\mathrm{TO}}(M,M))$, of
$${\mathcal C}^\infty(\hat{\mathcal A}_0(M)\times M\times \hat{\mathcal
B}(M),{\mathrm T}^r_s M).$$

So there remains the question: Why do we introduce slot 3 and how do transport
operators enter the scene? The answer is twofold:
%

%\pagebreak
\begin{itemize}
\item
Because it {\it works} (in German, we say ``Der Zweck
heiligt die Mittel'', i.e. ''The end justifies [sanctifies, literally] the
means'' in situations
like this), i.e.\ the resulting space $\hat{\mathcal E}^r_s(M)$ permits
sensible definitions of
\begin{itemize}
\item induced actions of diffeomorphisms $\mu:M\to N$,
\item natural extensions of Lie derivatives ${\mathrm L}_X$,
\item moderate and negligible elements and, finally
\item a space of generalized tensor fields $\hat{\mathcal G}^r_s(M)$ having all
the desired properties.
\end{itemize}
\item
The introduction of the $A$-slot for tensors is highly plausible---which the
remaining part of this article is devoted to convince the reader of.
\end{itemize}

Let us begin by reviewing the scalar case of embedding a (regular) distribution
given by a continuous function $g$ on $M$ into the basic space
$\hat{\mathcal E}(M)$, using formula (\ref{iota}) for the embedding $\iota$:
Pick $g\in {\mathcal C}(M)\subseteq{\mathcal D}'(M)$ and think of some $n$-form
$\omega$ which approximates the Dirac measure $\delta_p$ around $p\in M$; in
sloppy notation, $\omega(q)\approx\delta_p(q)$ for $p,q\in M$. Then
$$(\iota g)(\omega,p) = \int_M g(q)\omega(q)$$
{\it collects} values of $g$ around $p$ and forms a smooth {\it average} (note
that $\int\omega=1$!) as value for $(\iota g)(\omega,p)$. Here, $q\mapsto g(q)$
is a {\it scalar} valued function on $M$.

Now, if $g$ takes {\it tensors} of type $(r,s)$ as values, i.e., if $g$ is a
continuous section of ${\mathrm T}^r_s M$,
$$q\mapsto g(q)\in ({\mathrm T}^r_s)_q M,$$
then the $g(q)$'s do not live in the same linear space for different $q$!

In order to average them around $p$, we first have to ''gather'' them in $p$,
i.e. to shift each $g(q)$ from $({\mathrm T}^r_s)_p M$ to $({\mathrm
T}^r_s)_q M$.
This is accomplished by $A$ in the following way: For
$$g(q) = w_1(q)\otimes\dots\otimes w_r(q)\otimes
         \beta^1(q)\otimes\dots\otimes\beta^s(q)\ \in\ ({\mathrm
T}^r_s)_q M$$
(where $w_i(q)\in{\mathrm T}_qM$, $\beta^j(q)\in{\mathrm T}^*_qM$ for
$i=1,\dots,r$, $j=1,\dots,s$) we set
$$B^r_s(q,p)(g(q)):=A(q,p)w_1(q)\otimes\dots\otimes
               (A(p,q))^{\mathrm{ad}}\beta^s(q)\ \in\ ({\mathrm
T}^r_s)_p M.$$
(The notation $A^r_s(p,q)$ has to be saved for later use.) So we may form
\begin{equation}\label{neuesintegral}
(\iota^r_s g)(\omega,p,A) :=
     \int_M B^r_s(q,p)(g(q))\omega(q) \ \in\ ({\mathrm T}^r_s)_p M.
\end{equation}

In what follows, we will again simply write $\iota g$ for $\iota^r_s
g$.
Let us check the status of the objects in the above integrand carefully:
\begin{itemize}
\item
$B^r_s(q,p)(g(q))$ is an $(r,s)$-tensor at $p$, depending (smoothly) on
$q$;
\item
$\omega(q)$ (which $q$ viewed as variable) is a compactly supported $n$-form on
$M$ with unit integral.
\end{itemize}
So it seems that the integral on the right hand side of (\ref{neuesintegral})
is one of a ``new'' type (of course, only modulo the previous knowledge of
the reader), yet it is perfectly well-defined---just write it out in a chart in
the obvious way and check compatibility with chart changes.

As one can show, $\iota g$ as defined above depends smoothly on
$\omega,p,A$. (In fact, the proof of this statement represents one of the
technically most demanding parts of the forthcoming paper [3].) Thus for each
fixed pair $(\omega,A)$ we have that
$$(\iota g)(\omega,A):=[p\mapsto (\iota g)(\omega,p,A)]$$
defines
%a member of ${\mathcal T}^r_s M$, that is, by definition, it is
a smooth
tensor field of type $(r,s)$ on $M$, due to $(\iota g)(\omega,p,A)\in({\mathrm
T}^r_s)_p M$.

This strongly suggests the following choice for $\hat{\mathcal E}^r_s(M)$:
$$\hat{\mathcal E}^r_s(M):=\{
                 t\in{\mathcal C}^\infty(\hat{\mathcal A}_0(M)\times M\times
                 \hat{\mathcal B}(M),{\mathrm T}^r_s M)\,\mid\,
                 t(\omega,p,A)\in({\mathrm T}^r_s)_p M\}.$$
In particular, $p\mapsto t(\omega,p,A)$ is a member of ${\mathcal T}^r_s (M)$
for any fixed $\omega,A$.  As to the
inevitability of requiring smoothness in all three variables for the members of
the basic space,
see the remarks following formula (2) in [4].

Now, finally, we are going to pass from embedding {\it continuous}
$g$'s to embedding {\it distributional} tensor fields
$u\ \in\ {{\mathcal D}'}^r_s(M)$
into $\hat{\mathcal E}^r_s(M)$.

By definition of ${{\mathcal D}'}^r_s(M)$, $u$ takes (finite sums of)
tensors $\tilde t\otimes \omega$ ($\tilde t\in {\mathcal T}^s_r(M)$,
$\omega\in \Omega^n_c(M)$) as arguments.

Now what we need is a good formula for $(\iota u)(\omega,p,A)$. For a definition
of $\iota u$ in terms of $u$ we require something that $u$ can properly act
upon. We already have $\omega\in\Omega^n_c(M)$ from slot 1, so we still to have
to make some $\tilde t\in{\mathcal T}^s_r(M)$ enter the scene.

Fortunately, any $t\in{\mathcal T}^r_s(M)$ ($t=(\iota g)(\omega,A)$ in the case
at hand) is completely determined by specifying all contractions $t\cdot\tilde
t\in{\mathcal C}^\infty(M)$ where $\tilde t$ runs through ${\mathcal
T}^s_r(M)$. Hence we consider $(\iota g)(\omega,A)\cdot\tilde t$ defined
pointwise by

\begin{eqnarray*}
((\iota g)(\omega,A)\cdot\tilde t)(p)&=&
           (\iota g)(\omega,p,A)\cdot\tilde t(p)\\
           &=&\int_M B^r_s(q,p)(g(q))\cdot\tilde t(p)\,\,\omega(q)\\
           &=&\int_M g(q)\cdot (B^r_s(q,p))^{\mathrm{ad}}(\tilde
                     t(p))\,\,\omega(q)\\
           &=&\langle \underbrace{g(\,.\,)}_{\in{{\mathcal D}'}^r_s(M)}
                 ,\ \underbrace{(B^r_s(\,.\,,p))^{\mathrm{ad}}(\tilde
                     t(p))}_{\in{\mathcal T}^s_r(M)\mathrm{\ for\ fixed\ }p}
                  \ \otimes\underbrace{\omega(\,.\,)}_{\in \Omega^n_c(M)}\rangle
\end{eqnarray*}
(Note that in the third expression of the above calculation, 
$B^r_s(q,p)(g(q))$ and $\tilde t(p)$ are tensors of types $(r,s)$ and $(s,r)$,
respectively, hence their contraction is a scalar resp.\ a smooth function on
$M$. Therefore, the integrals above are usual integrals over $n$-forms rather
than of the ``new'' type discussed above.)

In the last expression above, we are now free to replace the regular
distribution $g$ by any $u\ \in\ {{\mathcal D}'}^r_s(M)$. This leads to our
definition of $\iota$, finally: Denoting
$(B^r_s(q,p))^{\mathrm{ad}}:({\mathrm T}^s_r)_p M\to({\mathrm T}^s_r)_q M$ by
$A^r_s(p,q)$, %[now it's happening!!]
we are led to define
\begin{eqnarray*}
(\iota u)(\omega,p,A)\cdot\tilde t(p)&:=&
                     ((\iota u)(\omega,A)\cdot\tilde t)(p)\\
            &:=&\langle\ u\ ,A^r_s(p,\,.\,)(\tilde
                       t(p))\,\otimes\,\omega(.)\ \rangle
\end{eqnarray*}
to obtain the desired embedding
$\iota=\iota^r_s:{{\mathcal D}'}^r_s(M)\to\hat{\mathcal E}^r_s(M)$.

With this definition of $\iota^r_s$, the requirement
$\hat\mu\circ\iota^r_s=\iota^r_s\circ\mu^*$ (for a diffeomorphism $\mu:M\to
N$) leads to a sensible definition of
$\hat\mu:\hat{\mathcal E}^r_s(N)\to\hat{\mathcal E}^r_s(M)$. The latter, in
turn, induces
$\hat{\mathrm{L}}_X:\hat{\mathcal E}^r_s(M)\to\hat{\mathcal
E}^r_s(M)$
satisfying $\hat{\mathrm{L}}_X\circ\iota^r_s=\iota^r_s\circ\mathrm{L}_X$
and $\hat{\mathrm{L}}_X\circ\hat\mu=\hat\mu\circ\mathrm{L}_{\mu_*X}$
(compare
[4]
%\cite{steingf07}
for the scalar case).

Corresponding to the above form of $\iota^r_s$ extending formula (\ref{iota})
to the tensor case, we also have the (much simpler) analog of (\ref{sigma})
for embedding smooth tensor fields into the new basic space:
$$ \sigma^r_s(f)(\omega,p,A):=f(p)\qquad(f\in{\mathcal T}^r_s(M)).$$
Preserving the product of smooth functions in the present context amounts to
preserving the tensor product of smooth tensor fields on $M$ when passing to
generalized functions via the embedding under discussion. This crucial goal of
the Colombeau approach again is achieved by an appropriate quotient construction
of moderate modulo negligible members of the basic space.

For a detailed account of the preceding introductory presentation, as well as
for an elaboration of
the
following concluding statement, we refer to the forthcoming paper
[3].
%\cite{gksv2}.

With the test for moderateness and negligibility from the scalar case
([2])
%(\cite{gksv})
suitably adapted to cope with slot 3, we finally arrive at 
$(\hat{\mathcal E}^r_s)_m(M)$,
$\hat{\mathcal N}^r_s(M)$ and
$$
\iota^r_s:{{\mathcal D}'}^r_s(M)\hookrightarrow
\hat{\mathcal G}^r_s(M):=(\hat{\mathcal E}^r_s(M))_m\,/\,\hat{\mathcal
             N}^r_s(M),$$
together with appropriate actions of diffeomorphisms and Lie derivaties on
$\hat{\mathcal G}^r_s(M)$ which naturally extend the corresponding notions on
${{\mathcal D}'}^r_s(M)$.

\vskip6pt
{\bf Acknowledgments.} The author wishes to express his gratitude to
the
organizers of the 
GF 07 conference in B\c{e}dlewo, in particular to Swiet\l ana
Minczewa-Kami\'nska and to Andrzej Kami\'nski.


\begin{thebibliography}{10}

\bibitem{gkos}
{Grosser, M., Kunzinger, M., Oberguggenberger, M., Steinbauer, R.,}
\newblock {\em Geometric Theory of Generalized Functions}. 
 Mathematics and its Applications 537.
\newblock Kluwer Academic Publishers, Dordrecht, 2001.

\bibitem{gksv}
{Grosser, M., Kunzinger, M., Steinbauer, R., Vickers, J. A.,}
\newblock A global theory of algebras of generalized functions.
\newblock {\em Adv. Math.} {\bf 166} (2002) 179--206.


\bibitem{gksv2}
{Grosser, M., Kunzinger, M., Steinbauer, R., Vickers, J. A.,}
\newblock A global theory of algebras of generalized functions II: tensor
distributions.
\newblock In preparation.


\bibitem{steingf07}
{Steinbauer, R.,}
\newblock A geometric approach to full Colombeau algebras.
\newblock In Kaminski, A., Pilipovi\'c, S., Oberguggenberger, M.,
editors, {\it Proceedings of the International Conference on
Generalized Functions 2007, Bedlewo, Poland}, Banach Center
Publications, to appear. arXiv:0710.2096v2 [math.FA] 

\bibitem{vw}
{Vickers, J. A., Wilson, J. P.,}
\newblock A nonlinear theory of tensor distributions.
\newblock {\em ESI-Preprint} (available electronically at
http://www.esi.ac.at./ESI-Preprints.html) {\bf 566}, 1998.

\end{thebibliography}
\end{document}